\documentclass[12pt]{amsart}

\usepackage[lite]{amsrefs}

\usepackage{amsmath}
\usepackage{amssymb}
\newtheorem{theorem}{Theorem}

\newtheorem{fact}[theorem]{Fact}

\newtheorem{thm}[theorem]{Theorem}

\theoremstyle{definition}

\theoremstyle{remark}

\numberwithin{equation}{theorem}

\newcommand{\reals}{{\mathbb R}}

\newcommand{\la}{\langle}
\newcommand{\ra}{\rangle}

\newcommand{\CM}{{\mathcal M}}

\newcommand{\complex}{\mathbb C}
\newcommand{\sub}{\subseteq}
\newcommand{\pf}{\noindent{\bf Proof }}

\newlength{\intlength}

\title[Not everything is standard]{A question of van den Dries and a
theorem of Lipshitz and Robinson;\\ Not everything is standard}

\author[Hrushovski]{Ehud~ Hrushovski}
\address{Department of Mathematics, Hebrew University,
Jerusalem, ISRAEL} \email{ehud@math.huji.ac.il}

\author[Peterzil]{Ya'acov~Peterzil}
\address{Department of Mathematics,
University of Haifa, Haifa, ISRAEL} \email{kobi@math.haifa.ac.il}

\begin{document}

\begin{abstract}
We use a new construction of an o-minimal structure, due to Lipshitz
and Robinson, to answer a question of van den Dries regarding the
relationship between arbitrary o-minimal expansions of real closed
fields and structures over the real numbers. We write a first order
sentence which is true in the Lipshitz-Robinson structure but fails
in any possible interpretation over the field of real numbers.
\end{abstract}

\maketitle
\date{}

An {\em o-minimal structure} is by definition an expansion $\mathcal
M$ of a linear ordering, such that every definable subset of the
linear ordering is a finite union of intervals whose end points are
in $\mathcal M\cup \{{\pm\infty}\}$.

Although o-minimal expansions of discrete linear orderings do exist
(e.g. $\la \mathbb Z,<,z\mapsto z+1\ra$), these were recognized
early on to have a relatively poor structure and
therefore, in the above definition, one often assumes that the
linear ordering is dense without endpoints.

As was shown in \cite{PS0}, an o-minimal structures is, at least
locally, one of the following three possibilities: It is degenerate
(basically an expansion of the linear ordering by unary functions),
an (interval in an) ordered vector space over an ordered division
ring, or an expansion of a real closed field.

Since ordered vector spaces over noncommutative ordered division
rings essentially cannot be further expanded while preserving
o-minimality (see \cite{PS0}), it is within the third possibility,
of expansions of real closed fields, where new o-minimal structures
can be found. Indeed, construction of new o-minimal structures is
usually carried out in this context, and much work has been done in
this direction in the past twenty years. However, until very
recently, all of these new o-minimal structures were expansions of
the field of real numbers $\la \reals,<,+,\cdot\ra$.

 Of course, we
know from model theory that every such structure $\mathcal M$ over
$\reals$ gives rise to an o-minimal expansion $\mathcal M^*$
of a nonstandard real
closed field. Now let $\mathcal
M_1^{*}$ be a structure with the same universe as $M^*$, and whose relations
are some of those definable with parameters in $M^*$.    $\mathcal
M_1^{*}$ may not be
elementarily equivalent to any structure  over the real numbers, for the simple reason
that it may  contain new constants  for ``infinitely large''
elements.

 A question arises. How can we
recognize that the theory of $\mathcal M_1^{*}$ has begun its model
theoretic life as a theory of a structure over
 the real numbers? And in general, is
it possible that {\em every}  o-minimal expansion of
a real closed field arises in this way?

Here is one possibility,
suggested by L. van den Dries in
\cite{VDD},  for a precise formulation to the above question:

Let $\phi(R_1,\ldots,R_n, f_1,\ldots,f_k)$ be a sentence in a
language $\mathcal L$ expanding the language of ordered rings, with
$R_1,\ldots,R_n, f_1,\dots,f_k$ all relation and function symbols
that are different than $\{<,+,\cdot\}$.
\vspace{3mm}

(*) \hspace{3mm} {\em Assume that $\mathcal N$ is
an o-minimal expansion of a real closed field in the language $\mathcal L$
such that
$\mathcal N\models \phi(R_1,\ldots,R_n,f_1,\ldots,f_k)$.  Does it
follow that:
$$\la \reals,<,+,\cdot\ra\models \exists
R_1,\ldots,R_n\,\exists f_1,\ldots,f_k \,\phi(R_1,\ldots,R_n,f_1,\ldots,f_k)?$$}

Notice that when  $\mathcal N$ equals $\mathcal M_1^{*}$ from the
above example then indeed the answer is positive.

For a given $\phi$ let us  denote by (**) the contra-positive to the
above. Namely the following transfer question:
\vspace{3mm}

\noindent(**)\hspace{3mm} {\em  If
$\phi(R_1,\ldots,R_n,f_1,\ldots,f_k)$ is true in every possible
expansion of $\la \reals,<,+,\cdot\ra$ then is it necessarily true
in every o-minimal expansion of a real closed field (in the language
$\mathcal L$)?} \vspace{2mm}

\noindent{\em Remark } As van den Dries points out, one gets a
negative answer to the above questions, if  the real closed field
assumption is omitted. Indeed, let $\phi$ be the statement: ``If
$\lambda_1\lambda_2$ are two continuous automorphisms of
$\la\reals,+\ra$ then $\lambda_1\lambda_2=\lambda_2\lambda_1$''.
This is true in every interpretation of $\lambda_1,\lambda_2$ over
the ordered group of the reals but fails in ordered vector spaces
over noncommutative ordered division rings. \vspace{2mm}

Here are some instances of sentences $\phi$ for which (**)
was (nontrivially) shown to have a positive answer:

1. (Invariance of Domain): If $f$ is a continuous injective function
from an open set in $\reals^n$ into $\reals^n$ then
it is an open map (\cite{WO}).

2.  If $G$ is a closed and bounded subset of $R^n$ and $\la
G,\star\ra$ is a topological group then $G$ has a torsion point
(see \cite{EO} for a precise count of the torsion points).

3. If $F$ is a function from the open unit disc in $\mathbb C$ into
$\mathbb C$ which is differentiable (with respect to $\mathbb C$)
then its derivative $F'(z)$ is differentiable as well (\cite{PS}).

It has already been pointed out that a positive answer to (*) would
yield some strong consequences. E.g., in \cite{BS}, Berarducci and
Servi showed that the decidability of the real exponential field
would follow. It would also have trivialized a significant part of
the theory of groups  in o-minimal structures. On the other hand, a
negative answer required a new procedure for constructing o-minimal
expansions of real closed fields different than the real numbers and
until recently none was known.

This has changed with the work of Lipshitz and Robinson (\cite{LR}),
where they constructed the following new o-minimal expansion of a
real closed field:

The underlying real closed field $R$ is the field of Puiseux series
in $t$ over $\reals$,
$$R=\bigcup_n\reals((t^{1/n})),$$ with $t$ an infinitesimal in the
ordering of the field.  If $p(\xi_1,\ldots,\xi_n)$ is
a formal power series over $\reals$ (the ring of all
such power series is denoted by  $\reals[[\xi_1,\ldots,\xi_n]]$) then
$p(\xi_1,\ldots,\xi_n)$ converges on the infinitesimal cube
$[-t,t]^n\sub R^n$. Denote by  $f_p(x_1,\ldots,x_n)$  the
corresponding function, which is set to be zero outside the box
$[-t,t]^n$. The theorem of Lipshitz and Robinson says:
\begin{thm} The
structure $\mathcal M=\la R,f_p\ra_{p\in
\reals[[\xi_1,\ldots,\xi_n]]}$ is o-minimal.
\end{thm}

Here are some basic properties:
\begin{fact}\label{fact1}
\begin{enumerate}
\item The map $p\mapsto f_p$ is an embedding of rings from the ring
 $\reals[[\xi_1,\ldots,\xi_n]]$ into the
ring of definable functions on $[-t,t]^n$.
\item For $p\in \reals[[\xi_1,\ldots,\xi_n]]$, and $q=\partial
p/\partial x_i$ the formal derivative with respect to $\xi_i$, we
have $\partial f_p/\partial x_i=f_q$ on $-[t,t]^n$, where the
partial derivative of $f_p$ is taken with respect to the real closed field $R$.
\end{enumerate}
\end{fact}
\pf (1) is standard. For (2), consider for simplicity the 1-variable
case. Then, it follows from (1), that
 for small $x$ and $h\neq 0$,
$$\frac{f_p(x+h)-f_p(x)}{h}=f_{\frac{p(\xi+h)-p(\xi)}{h}}(x).$$
The result  easily follows.\qed

Denote by $K=R(\sqrt{-1})$ the algebraic closure of $R$, identified
as usual with $R^2$. Notice that the product topology induced on $K$
from $R$ is the same as the valuation topology when we identify $K$
with the field of Puiseux series over $\mathbb C$. We denote by
$D\sub K$ the open disc around of radius $t$ (inthe sense of $R^2$)
around $0\in K$. As in \cite{PS}, we say that an $\CM$-definable
function $F$ from an open set $U\sub K$ into $K$ is {\em
$K$-differentiable at $z_0\in U$} if the limit
$$\lim_{z\to 0}\frac{F(z_0+z)-F(z_0)}{z}$$ exists in $K$.
We call this limit $F'(z_0)$. Notice that the definition through limits
gives a formula $\psi_{der}(F,z_0)$
in the language of ordered rings, augmented
by a symbol for $F$ (or more precisely by two functions symbols for the real
and imaginary parts of $F$), such that for any structure
$\la \reals,<,+,\cdot, F\ra$ and any $z_0\in \mathbb C$, the function $F(z)$
is complex differentiable at $z_0$ if and only if $\psi_{der}(F,z_0)$ holds.

We now return to our field of Puiseux series $R$ and its algebraic
closure $K$. If $p(\zeta)$ is  a formal power series in
$\reals[[\zeta]]$ then $p(z)$ converges, in the topology induced on
$K$ from $R$,  for all $z\in D\sub K$. Let $F_p(z)$ be the
corresponding function from $D$ into $K$.  Notice that $F_p$
 as two coordinate functions $(f_{p_1},f_{p_2})(x,y)$, where
$p_1, p_2$ are themselves power series in $\reals[[\xi_1,\xi_2]]$
(indeed, this follows from the fact that each map $z\mapsto z^n$ can
be written,  in $R$-coordinates, as $(q_1(x,y),q_2(x,y))$ where
$q_1$ and $q_2$ are homogeneous polynomials  of degree $n$).
Therefore, the function $F_p(z)$ is definable in $\mathcal M$.

\begin{fact}\label{fact2}
\begin{enumerate}

\item The map $p\mapsto f_p$ is an embedding of rings from the
ring $\reals[[\zeta]]$ into the ring of definable functions from $D$
into $K$.
\item For $p\in \reals[[\zeta]]$, the function $F_p(z)$ is
$K$-differentiable. If $p'(\zeta)$ is the formal derivative with
respect to $\zeta$, we have $F_p'(z)=F_{p'}(z)$ for all $z$ in $D$.
\end{enumerate}
\end{fact}
\pf The proof of (1) follows from Clause (1) in Fact \ref{fact1}.
The proof of (2) is just like the
proof of Clause (2) above.\qed

 We can now produce a sentence $\phi$ for which the answer to (*) is negative:
The signature contains two 2-ary function symbols $f_1$ and $f_2$.

Let
$\phi(f_1,f_2)$ be the following sentence (using $F(z)$
for simplification instead of
$(f_1,f_2)(x,y)$):

{\em There exists $r>0$, such that for all $z=x+iy$, if $|z|<r$ then
$\psi_{der}(F,z)$ and
\begin{equation} \label{eq1} F(z)=z^2F'(z)+z. \end{equation}}

\begin{fact}
\begin{enumerate}
\item Let $p(\zeta)=\Sigma_{n=1}^{\infty}(n-1)!{\zeta}^n.$ Then $F_p(z)$
is a solution to \ref{eq1} on some open neighborhood of $0\in K$.
\item The sentence $\phi(f_1,f_2)$ is false in
$\la\reals,<,+,\cdot, f_1,f_2\ra$
for every possible interpretation of $f_1,f_2$.
\end{enumerate}
\end{fact}
\pf (1) It is easy to verify that $p(\zeta)$ is a {\em formal}
solution to (\ref{eq1}).
It follows from Fact \ref{fact2} that $F_p$ is a solution to
the same equation on $D$.

For (2), notice that if $F(z)$ were a complex differentiable
function which solves \ref{eq1} then its Taylor series at $0$ must be
$p(\zeta)$. However, this power series is divergent at every
nonzero $z\in \complex$.\qed

The sentence $\phi(f_1,f_2)$ is therefore true in $\mathcal M$, with
$f_1^{\mathcal M},f_2^{\mathcal M}$ interpreted as the real and
imaginary parts of $F_p$ (for $p$ as in the last Fact), but fails
over the real numbers in every possible interpretation of $f_1,f_2$.
\vspace{3mm}


\noindent{\bf A question } Because of its strong potential
consequences we are tempted to re-formulate van den Dries' original
question as follows:

Find a restricted class $\mathcal K$ of o-minimal expansions of real
closed fields such that every finite theory $T$ in a language
$\mathcal L$ expanding the langauge of real closed fields, which
holds in some o-minimal expansion of a real closed field, has an
interpretation in one of the structures in $\mathcal K$.

Another possible variation on the requirement from the above
$\mathcal K$ is the following local version:

Consider a sentence $\phi$ in a language $\mathcal L$ expanding real
closed fields, and an o-minimal $\mathcal L$-structure $\CM$.  For
every $t>0$ in $\mathcal M$, consider the sentence $\phi_t$ obtained
from $\phi$ by restricting all quantifiers and relations to
cartesian products of $[-t,t]$. By o-minimality, the truth value of
$\phi_t$ in $\mathcal M$ stabilizes as as $t$ approaches $0$. We
ask: Given $\phi$ and $\mathcal M$, is there a structure ${\mathcal
M}_{\phi}$ in $\mathcal K$ and an interpretation of all the symbols
in $\phi$ in ${\mathcal M}_{\phi}$ such that the limit truth value
of $\phi_t$ is the same in $\mathcal M$ and in $\mathcal M_{\phi}$?

We feel that the Lipshitz-Robinson model described above could play
a significant role in finding an appropriate $\mathcal K$.

\vspace{3mm}

We first realized that an equation such as \ref{eq1} should exist by
considering Sofia Kovalevskaya's example of an analytic differential
equation whose solutions are not analytic.  The actual example in
this paper was arrived at after a series of simplifications and
corrections.

\end{document}